\documentclass[11pt,a4paper]{article}
\usepackage{amsmath,amssymb}

\usepackage{amsmath}

\usepackage{color}
\usepackage[colorlinks,linkcolor=blue,citecolor=red]{hyperref}

\DeclareMathSymbol{\bbbr}{\mathalpha}{AMSb}{"52}
\DeclareMathSymbol{\bbbc}{\mathalpha}{AMSb}{"52}

\newtheorem{theorem}{Theorem}

\newtheorem{corollary}[theorem]{Corollary}

\newtheorem{definition}[theorem]{Definition}

\newtheorem{lemma}[theorem]{Lemma}

\newtheorem{proposition}[theorem]{Proposition}

\begin{document}

\title{Integrable geodesic flow in 3D and webs of maximal rank}

\author{{\Large Sergey I. Agafonov}\\
\\
Department of Mathematics,\\
S\~ao Paulo State University-UNESP,\\ S\~ao Jos\'e do Rio Preto, Brazil\\
e-mail: {\tt sergey.agafonov@gmail.com} }
\date{}
\maketitle
\unitlength=1mm

\vspace{1cm}

\begin{abstract} We characterize  geodesic flows, admitting two commuting quadratic integrals with common principal directions, in terms of the geodesic 4-webs such that the tangents to the web leaves are common zero directions of the integrals. We prove that, under some natural geometric hypothesis, the metric is of Stäckel type.
\bigskip

\noindent MSC: 53A60, 37J06.

\bigskip

\noindent
{\bf Keywords}: web rank, integrable geodesic flow, Stäckel metric.
\end{abstract}



\section{Introduction}

Integrable geodesic flow is a classical theme of analytic mechanics and  differential geometry. Its foundation was laid  in the 19th century by works of Jacobi \cite{J-84}, Darboux \cite{D-91}, Dini \cite{D-69},  Koenigs \cite{K-96}, Stäckel \cite{S-93}, and others. The most developed is the theory of  integrals  polynomial in momenta, also called {\it Killing tensors}. Even in the local setting, this problem is well motivated  by the following observations: a) if there is a locally defined  analytic integral then each  homogeneous in momenta part of its Taylor expansion is an integral  polynomial in momenta, b) the submanifold of zero velocities  in the tangent bundle  is the set of singular points of the geodesic flow, the existence of regular integrals at singular point being  nontrivial.

The geodesic flow is Hamiltonian, therefore, by Liouville's Theorem,   it is enough to find one first integral independent of the Hamiltonian to integrate the flow on a surface. For three-dimensional Riemannian manifold, one needs two commuting integrals such that they and the Hamiltonian form a functionally independent triple.

Recent interest to integrable geodesic flow on surfaces is due to a new discovered relation to infinite-dimensional integrable systems of partial differential equations (PDE). Namely, if one considers the metric and the integral  as unknown then the system of PDE on metric and integral is of hydrodynamic type and is linearizable by Tsarev's {\it generalized hodograph method}: it is diagonalizable and it possesses infinitely many conservation laws  (see \cite{K-82,T-97,BM-11,MP-17,T-85, T-90}).

Integrable geodesic flow in higher dimensions were intensively studied for invariant metrics on Lie groups, the topic also being classical since the studies of rigid body motion by Euler. Another known example is the so called Stäckel metric. Importance of this example is due to its various non-trivial relations to integrable differential equations.

Thus, Stäckel himself \cite{S-91,S-93} used these diagonal metrics to find complete integral of Hamilton-Jacobi equation by separation of variables. Recall that the complete integral allows to write finite (i.e. not differential) equations of geodesics in terms of differentiating and inverting of maps.    He also indicated the set of commuting quadratic integrals of the corresponding geodesic flow.

Eisenhart \cite{E-34} proved a converse of Stäckel's result: if the geodesic flow on an $n$-dimensional Riemannian manifold admits $n-1$ commuting quadratic integrals $I_k$ such that 1) these integrals, written as functions of velocities, and the metric $g$ can be simultaneously diagonalized in  some orthogonal coordinate system and 2) all $n-1$ characteristic equations $\det (I_k-\lambda g)=0$ have simple roots  then the metric in these coordinates takes the Stäckel form.

Benenti \cite{B-92} showed how to recover $n-1$ commuting quadratic integrals from a Killing tensor $L$  of {\it gradient type} with vanishing Nijenhuis torson. (The tensor $L$ is of gradient type if the Poisson bracket $\{L,H\}$ factors into $H$ and a form linear in momenta. For the Benenti construction to work, this linear factor must be of special form.) The obtained integrals satisfy the condition of the Eisenhart theorem and therefore the metric is the one of Stäckel. These results allowed Benenti to extend the separation of variables approach to Hamiltonians with potential terms \cite{Bo-92,B-97}.

Ibort, Magri, and  Marmoin \cite{IMM-00} interpreted the Benenti results in terms of the Gel'fand-Zakharevich bi-Hamiltonian structure on the extension of the phase space for the corresponding Hamiltonian system. Ferapontov and Fordy \cite{FF-97} found a remarkable relation of Stäckel metrics to integrable linearly degenerate systems of hydrodynamic type.

The objective of this paper is to contribute to the reach structure behind the Stäckel metrics, namely, we interpret the existence of Stäckel metric in terms of the web theory. For cubic, quadratic and fractional-linear integrals on surfaces,
 the geometry of corresponding webs was revealed in \cite{Ah-21,Aq-21,AA-23}. For polynomial integrals, the web is constructed as follows: we write the integral $I$ of degree $d$ as a homogeneous form in velocities and consider the equation  $I=0$. For each point on the surface, this equation gives $d$  directions. If these directions are real then their integral curves are geodesics. Thus we have $d$ foliations by geodesics, i.e. geodesic $d$-web. For cubic integrals, this web is {\it hexagonal} and conversely, existence of a hexagonal geodesic 3-web implies the existence of a cubic integral with real zero directions. This result \cite{Ah-21} solves the following Blaschke problem \cite{B-55}: to find an intrinsic criterion for existence of hexagonal geodesic 3-webs, thus  generalizing the Graf and Sauer Theorem \cite{GS-24} for metrics of non-constant gaussian curvature and providing a new "dynamical" proof of this Theorem \cite{Aq-21}. For quadratic integrals with real zero directions we get a geodesic 2-web (or geodesic net) $G_2$. A net does not have topological invariants. Using the metric structure, we define the {\it bisector net}, whose leaf directions bisects the angles between the leaves of $G_2$. Then each 3-subweb of the 4-web formed by $G_2$ and $B_2$ is hexagonal and the converse is also true (for more detail see \cite{Aq-21}).

This construction of geodesic webs via polynomial integrals of geodesic flows generalizes to any dimension. Consider two quadratic commuting integrals $I_2,I_3$, written as homogeneous forms in velocities, in 3-dimensional Riemannian manifold $M$. Suppose that for $\lambda,\mu \in \mathbb{R}$, the eqautions
$$
I_2-\lambda g=I_3-\mu g=0
$$
define 4 different real points in  $PT_mM$, thus giving 4 real direction fields $\tau_i$. Integral curves of these fields are geodesics, we obtain a geodesic 4-web. For Stäckel metric and integrals, the reflection in tangent planes to the coordinate surfaces maps any $\tau_i$ to the other three. Now fix one $\tau_i$  and consider the 4-web composed of integral curves of $\tau_i$ and of 3 coordinate surfaces. Then this 4-web is of maximal rank (see Proposition \ref{rank}).

Our first main result claims that the converse is also true.
\begin{theorem}\label{webstack}
Suppose that a three-dimensional Riemannian manifold $M$ carries a geodesic 4-web $w_g$ such that at each point $m\in M$ there are 3 mutually orthogonal planes $\pi_i$, $i=1,2,3,$ such that
\begin{enumerate}
\item the distributions $\pi_i$ are integrable,
\item the  reflections at $\pi_i$ map one direction of the geodesic web $w_g$ to the other three,
\item the rank of the 4-web $w_i$ formed by geodesics of one foliation of $w_g$ and by three foliations by integral manifolds of $\pi_1,\pi_2,\pi_3$ is maximal and equals to 2.
\end{enumerate}
Then the metric is of Stäckel type and $\pi_i$ are tangent to  the Stäckel coordinate surfaces.
\end{theorem}

Our second main result is relaxing of the hypothesis of Eisenhart's Theorem \cite{E-34} for the dimension three. Let $I$ be quadratic in velocities integral, then the roots $\lambda$ of the characteristic equation
$$
\det (I-\lambda g)=0
$$
are called the {\it principal invariants} and the vectors $v\in T_mM$, verifying
$$
 (I-\lambda g)v=0
$$ are called the {\it principal directions} of the corresponding symmetric tensor $I$ (see \cite{E-34} p. 110). Note that there is an orthogonal basis of the principal directions even for multiple principal invariants. The principal directions are called {\it normal} if their normal distributions are integrable. We show that one can omit normality of the principal directions assumed by the Eisenhart theorem.
\begin{theorem}\label{normal}
If the geodesic flow  on a three-dimensional manifold admits   three independent commuting  integrals, quadratic in momenta (including the Hamiltonian),  such that
 the integrals share the same principal directions $n_i$ then all
 the distributions, orthogonal to $n_i$, are integrable and the metric is of Stäckel type.

\end{theorem}

All the results of this paper are local, all functions, surfaces, fields, and other objects are smooth.
\section{Stäckel metrics and geodesic webs}
Stäckel metrics  on $n$-dimensional manifold $M$ with local coordinates $x_i$, $i=1,...,n$ are constructed from $n^2$ functions $\varphi_{ij}(x_i)$, each depending on one variable:
\begin{equation}\label{staeckelmetric}
g= \sum_i\frac{\Delta}{\Phi^{i1}}dx_i^2,
\end{equation}
where $\Delta:= \det \varphi_{ij}\ne 0$ and $\Phi^{ij}$ is the cofactor of $\varphi_{ij}$ in $\Delta$.
The geodesic flow of this metric admits $n-1$ quadratic integrals  in involution $I_k$, $k=2,...,n$ (see \cite{S=91,S-93} or \cite{E-49}), corresponding to the quadratic forms
$$
I_k= \sum_i\rho_{ik}\frac{\Delta}{\Phi^{i1}}dx_i^2,\ \ \ \mbox{where}\ \ \ \rho_{ik}=\frac{\Phi^{ik}}{\Phi^{i1}}.
$$
In what follows we call $x_i$ Stäckel coordinates.
For $n=3$ and $x_1=x$, $x_2=y$, $x_3=z$ let the matrix $\varphi_{ij}$ be
$$
\left(
\begin{array}{ccc}
E(x) & K(x) & P(x)\\
F(y) & L(y)  & Q(y)\\
G(z) & M(z) & R(z)
\end{array}
\right).
$$
Then $\Delta = ELR+FMP + GKQ - EMQ - FKR - GLP$ and
\begin{equation}
g=\frac{\Delta}{LR - MQ}dx^2 + \frac{\Delta}{MP-KR}dy^2 + \frac{\Delta}{KQ - LP}dz^2,
\end{equation}
\begin{equation}
I_2=\frac{\Delta (GQ-FR)}{(LR - MQ)^2}dx^2 + \frac{\Delta (ER - GP)}{(MP-KR)^2}dy^2 + \frac{\Delta (FP-EQ)}{(KQ - LP)^2}dz^2,
\end{equation}
\begin{equation}
I_3=\frac{\Delta (FM - GL)}{(LR-MQ)^2}dx^2 + \frac{\Delta (GK-EM)}{(MP-KR)^2}dy^2 + \frac{\Delta (EL - FK)}{(KQ-LP)^2}dz^2.
\end{equation}

  \medskip
\noindent {\bf Remark.} Adapting local coordinates, one can reduce the number of parameterizing functions to $n(n-1)$. For $n=3$ we define $a=\frac{E}{P}$, $k=\frac{K}{P}$, $b=\frac{F}{Q}$, $l=\frac{L}{Q}$, $c=\frac{G}{R}$, $m=\frac{M}{R}$.
Then $\Delta=PQR\cdot \delta$, where $\delta=a(l-m)+b(m-k)+c(k-l),$ and
$$
\frac{\Delta}{LR - MQ}=\frac{\delta P}{l-m}, \ \ \  \frac{\Delta}{MP-KR}=\frac{\delta Q}{m-k}, \ \ \ \frac{\Delta}{KQ - LP}=\frac{\delta R}{k-l}.
$$
Now in the  new local coordinates $\bar{x}$, $\bar{y}$, $\bar{z}$  chosen  so that
$$
P(x)dx^2=a(x)d\bar{x}^2, \ \ \ Q(y)dy^2=b(y)d\bar{y}^2, \ \ \ R(z)dz^2=c(z)d\bar{z}^2,
$$
the metric is
\begin{equation}\label{stform}
g=\frac{\delta a}{l - m}d\bar{x}^2 + \frac{\delta b}{m-k}d\bar{y}^2 + \frac{\delta c}{k - l}d\bar{z}^2.
\end{equation}
  \medskip
Let us fix  $\lambda$ and $\mu$ in a way that the two equations
\begin{equation}\label{webdirectioneq}
I_2-\lambda g=I_3-\mu g=0
\end{equation}
define 4 different real points in $PT_mM$. Then  these points define 4 real direction fields $\tau_i$, $i=1,2,3,4$ (at least  locally).
\begin{proposition}\label{mirrors}
Integral curves of $\tau_i$ are geodesic, at each point there are 3 mutually orthogonal planes $\pi_i\subset T_mM$, $i=1,2,3,$ such that reflections at $\pi_i$ map one direction $\tau_k$ to the other three. Moreover, the distributions $\pi_i$ are integrable.
\end{proposition}
{\it Proof:} The first claim follows from the uniqueness  of an integral curve passing through a point: the forms $I_2-\lambda g$, $I_3-\mu g$, evaluated on a vector field with the direction $\tau_k$, vanish   along the geodesic outgoing from the point at the direction  $\tau_k$.

The planes $\pi_i$ are tangent to the coordinate surfaces, symmetry properties follow from the diagonal form of $g,I_2,I_3$ in coordinates $x,y,z$.
\hfill $\Box$\\

\smallskip

The geodesic 4-web defined by (\ref{webdirectioneq}) has a more subtle property expressed in terms of the web rank.
\begin{definition}
Abelian relation of a d-web $w$ formed by foliations $\mathcal{F}_1,\mathcal{F}_2,...,\mathcal{F}_d$ is a d-tuple of one-forms $(\sigma_1,\sigma_2,...,\sigma_d)$ such that:
\begin{enumerate}
\item each  $\sigma_i$ vanishes on the leaves of $\mathcal{F}_i$,
\item $d\sigma_i=0, \ \ \ i=1,...,d,$
\item $\sigma_1+\sigma_2+...+\sigma_d=0.$
\end{enumerate}
The web rank $r_w$ is the real dimension of the real vector space $\mathcal{A}_w$ of abelian relations.
\end{definition}
Web rank can be infinite, it is defined also for webs whose leaves have different codimensions for different foliations. Consider the 4-web $w_i$ formed by coordinate surfaces $x=const$, $y=const$, $z=const$ and by integral curves of the field $\tau_i$.
\begin{proposition}\label{rank}
The rank of the web $w_i$ is at least two.
\end{proposition}
{\it Proof:}
Equations (\ref{webdirectioneq}) give
$$
\frac{(KQ - LP)^2}{E+\lambda K+\mu P}dx^2= \frac{(LR - MQ)^2}{G+\lambda M+\mu R}dz^2,\ \ \ \ \
\frac{(KQ - LP)^2}{F+\lambda L+\mu Q}dy^2=\frac{(KR - MP)^2}{G+\lambda M+\mu R}dz^2.
$$
Since the directions $\tau_i$ are real, the expressions $E+\lambda K+\mu P$, $F+\lambda L+\mu Q$ and  $G+\lambda M+\mu R$ have the same sign in some neighborhood of a point. Define closed forms
\begin{equation}
\begin{array}{c}
\omega_1=\frac{dx}{\sqrt{|E(x)+\lambda K(x)+\mu P(x)|}}, \\
\omega_2=\frac{dy}{\sqrt{|F(y)+\lambda L(y)+\mu Q(y)|}}, \\
\omega_3=\frac{dz}{\sqrt{|G(z)+\lambda M(z)+\mu R(z)|}}.
\end{array}
\end{equation}
Now equations (\ref{webdirectioneq}) write as
$$
(KQ-LP)\omega_1\pm (LR-MQ)\omega_3=0,\ \ \ (KQ-LP)\omega_2\pm (KR-MP)\omega_3=0,
$$
where all possible choices of signs give all $\tau_i$. The above two equations are equivalent to
\begin{equation}\label{directions}
\begin{array}{l}
e_k K(x)\omega_1+ e_l L(y)\omega_2+M(z)\omega_3=0,\\
\\
 e_k P(x)\omega_1+ e_l Q(y)\omega_2+R(z)\omega_3=0,
\end{array}
\end{equation}

where $e_k=\pm 1$, $e_l=\pm 1$.
Thus, for $e_k=e_l=1$, the 4-tuples
$$
a_1=(- K(x)\omega_1, -L(y)\omega_2, -M(z) \omega_3, K(x)\omega_1 + L(y)\omega_2+M(z)\omega_3 )
$$
$$
a_2=(- P(x)\omega_1, -Q(y)\omega_2, -R(z) \omega_3, P(x)\omega_1 + Q(y)\omega_2+R(z)\omega_3 )
$$
are independent abelian relations for the corresponding web $w_i$. For the other choices of signs, one easily adapts $a_1$ and $a_2$.
\hfill $\Box$\\
In the next section we will see that the rank 2 is maximal possible for the considered type of webs.

\section{Maximal rank webs}\label{secmaxrank}
\begin{theorem}
The maximal rank of a non-degenerate 4-web  in $\mathbb{R}^3$,  composed of three foliations by surfaces  and  one foliation by curves,  is two.
\end{theorem}
{\it Proof:} Suppose that the web rank is at least two. We can choose local coordinates $x,y,z$ in $\mathbb{R}^3$ so that  the web surfaces are given by equations $x=const$, $y=const$, $z=const$.   Moreover, without loss of generality, we can assume that the basic abelian relations are $a_1=(-dx,-dy,-dz,dx+dy+dz)$ and  $a_2=(-f(x)dx,-j(y)dy,-h(z)dz,f(x)dx+j(y)dy+h(z)dz)$. Then the vector field $v=(h-j)\partial_x+(f-h)\partial_y+(j-f)\partial_z$ is tangent to the foliation curves. Let $a=(-u(x)dx,-v(y)dy,-w(z)dz,u(x)dx+v(y)dy+w(z)dz)$ be an abelian relation, then
\begin{equation}\label{abeleq}
u(x)[h(z)-j(y)]+v(y)[f(x)-h(z)]+w(z)[j(y)-f(x)]=0.
\end{equation}
 Taking second order mixed derivatives  of the left hand side expression, one gets
$$
f'(x)v'(y)-j'(y)u'(x)=0, \ \ \ j'(y)w'(z)-h'(z)v'(y)=0,\ \ \ h'(z)u'(x)-f'(x)w'(z)=0.
$$
If at least one of the functions $f,j,h$ is not constant then, without loss of generality, one assumes $h'(z)\ne 0$ and gets
$$
u'(x)=\frac{w'(z)}{h'(z)}f'(x) \ \ \  \mbox{and} \ \ \ v'(y)=\frac{w'(z)}{h'(z)}j'(y),
$$
hence $\frac{w'(z)}{h'(z)}=k=const$ or $f'(x)=j'(y)=u'(x)=v'(y)=0.$

In the former case holds
\begin{equation}\label{uvw}
u(x)=kf(x)+n, \ \ \  v(y)=kj(y)+m, \ \ \  w(z)=kh(z)+l
\end{equation}
with constant $l,m,n$. Therefore $(-ndx,-mdy,-ldz,ndx+mdy+ldz)=a-ka_2$ is also an abelian relation. Substituting (\ref{uvw}) with $k=0$ into (\ref{abeleq}), differentiating by $z$ and taking into account $h'\ne 0$ we arrive at $m=n$.
If  the abelian relation $a-ka_2-na_1=(0,0,(n-l)dz,(l-n)dz)$ is not trivial then   the vector field $v$ belongs to the distribution $dz=0$ and the web is degenerate. Therefore $l=n$ and the space of abelian relation is two-dimensional.

In the latter case, the functions $f,j,u,v$ are constant. Differentiating (\ref{abeleq}) by $z$ and taking into account $h'\ne 0$ we get or $\frac{w'(z)}{h'(z)}=const$ or $f=j$, $u=v$ and web is again degenerate.

Finally, suppose that $f'=j'=h'=0$. Then the functions $f,j,h$ are constant and, since the web rank is not less than two, at least two differences of $f-j,j-h,h-f$ are not zero. Suppose that $h\ne j$ and $f\ne h$.
Differentiating  (\ref{abeleq}) by $x$ we get $u'(x)=0$, similarly $v'(u)=0$, and $u,v$ are also constant. If $w'=0$ then $w$ is also constant and equation (\ref{abeleq}) means that $a$ is a linear combination of $a_1,a_2$. Otherwise,
 differentiating  (\ref{abeleq}) by $z$ we get $f=j$. Then $a_2-fa_1=(0,0,(f-h(z))dz,(h(z)-f)dz)$ is an abelian relation and the web is degenerate.
\hfill $\Box$ \\

We will need also a rank bound for the  following  type of 3-webs.

\begin{lemma}
The maximal rank of a non-degenerate 3-web  in $\mathbb{R}^3$,  composed of two foliations by surfaces  and  one foliation by curves,  is one.
\end{lemma}
{\it Proof:} Suppose that the web rank is one. We can choose local coordinates $x,y,z$ in $\mathbb{R}^3$ so that  the web surfaces are given by equations $x=const$, $y=const$.   Without loss of generality, we can assume that the abelian relation is $a_1=(-dx,-dy,dx+dy)$. Then the vector field, tangent to the foliation curves, has the form $v=f\partial_x -f\partial_y+h\partial_z$ . Since the web is non-degenerate, $f\ne 0$ and we can assume $f=1$.
Let $a=(-u(x)dx,-v(y)dy,u(x)dx+v(y)dy)$ be an abelian relation.
Then $u(x)-v(y)=0$ hence $u=v=const$, $a=ua_1$ and the rank is one.
\hfill $\Box$\\

There is a criterion for the rank of 3-webs described in Lemma to have the  rank one. Let the 2D foliations be defined by the kernels of 1-forms $\bar{\omega}_1$ and $\bar{\omega}_2$ and the foliation by curves by the kernel intersection of $\bar{\omega}_3$ and $\bar{\omega}_4$. One can assume the forms $\bar{\omega}_1,\bar{\omega}_2,\bar{\omega}_3$ to constitute a basis and normalize, up to a common factor, the chosen forms to satisfy
\begin{equation}\label{norm221}
\bar{\omega}_4=\bar{\omega}_1+\bar{\omega}_2+\bar{\omega}_3.
\end{equation}
The forms verify the structure equations
\begin{equation}
\begin{array}{l}\label{structureweb}
         d\bar{\omega}_1= \bar{a} \bar{\omega}_1 \wedge\bar{\omega}_2 +\bar{ b} \bar{\omega}_1 \wedge\bar{ \omega}_3 ,\\

         d\bar{\omega}_2 = \bar{p} \bar{\omega}_2 \bar{\wedge \omega}_3 + \bar{q} \bar{\omega}_2 \wedge \bar{\omega}_1 ,\\

         d\bar{\omega}_3= \bar{r} \bar{\omega}_3 \wedge \bar{\omega}_1 + \bar{s} \bar{\omega}_3 \wedge \bar{\omega}_2 + \bar{M} \bar{\omega}_1 \wedge \bar{\omega}_2.
\end{array}
\end{equation}
Then the abelian relations have the form
\begin{equation}\label{intfac}
\bar{K}\bar{\omega}_1+\bar{K}\bar{\omega}_2+\bar{K}(\bar{\omega}_3-\bar{\omega}_4)=0,\ \mbox{with}\ d(\bar{K}\bar{\omega}_1)=d(\bar{K}\bar{\omega}_2)=0.
\end{equation}

\begin{lemma}\label{2s1c}
The rank of a 3-web  in $\mathbb{R}^3$,  composed of two foliations by surfaces  and  one foliation by curves,  is one if and only if
$d(\bar{\omega}_1\wedge \bar{\omega}_2)=0$ and $d(\gamma)=0$, where
$\gamma=\bar{q}\bar{\omega}_1+\bar{a}\bar{\omega}_2+\bar{b}\bar{\omega}_3$.
\end{lemma}
{\it Proof:} Differential of any function can be written in the basis $\bar{\omega}_1$, $\bar{\omega}_2$, $\bar{\omega}_3$, for example,
$$
d\bar{K}=\bar{K}_1\omega_1+\bar{K}_2\bar{\omega}_2+\bar{K}_3\bar{\omega}_3.
$$
The common integrating factor $\bar{K}$ defined by (\ref{intfac}) exists if and only if $\frac{\bar{K}_3}{\bar{K}}=\bar{b}=\bar{p}$, $\frac{\bar{K}_1}{\bar{K}}=\bar{q}$, $\frac{\bar{K}_2}{\bar{K}}=\bar{a}$ and the form $\gamma$ is closed. Hence the claim.
\hfill $\Box$

\begin{definition}
The form $\gamma$ is called the connection form of a web 3-web  in $\mathbb{R}^3$,  composed of two foliations by surfaces  and  one foliation by curves if $$
d\bar{\omega}_1=\gamma\wedge \bar{\omega}_1, \ \ \ d\bar{\omega}_2=\gamma\wedge \bar{\omega}_2.
$$
The curvature of this 3-web is  $d(\gamma)$.
\end{definition}
The definition is motivated by the connection of a planer 3-web. Note that\\ 1) the connection of the considered web type exists if and only if $d(\bar{\omega}_1\wedge \bar{\omega}_2)=0$,\\ 2) this condition is invariant by normalizations for integrable distributions $\bar{\omega}_1=\bar{\omega}_2=0,$\\ 3) the connection depends on normalization (\ref{norm221}) whereas the curvature does not.
\begin{corollary}
The rank of a 3-web  in $\mathbb{R}^3$,  composed of two foliations by surfaces  and  one foliation by curves,  is one if and only if
the web curvature vanishes.
\end{corollary}

\section{Stäckel metrics from  webs of maximal rank}
In this section we prove the converse claim of Proposition \ref{rank}. \\
{\bf Proof of Theorem \ref{webstack}:} Since the distributions $\pi_i$ are integrable, there are functions whose level surfaces are integral manifolds of these distributions. Moreover, one can choose such functions as local coordinates $x,y,z$ so that $dx|_{\pi_1}=0$,  $dy|_{\pi_2}=0$, $dz|_{\pi_3}=0$ and the basic abelian relation of $w_i$ are $a_1=(-dx,-dy,-dz,dx+dy+dz)$ and  $a_2=(-k(x)dx,-l(y)dy,-m(z)dz,k(x)dx+l(y)dy+m(z)dz)$. Then  the vector field
$$
v=\xi\partial_x+\eta\partial_y+\zeta\partial_z=(l-m)\partial_x+(m-k)\partial_y+(k-l)\partial_z
$$
 is tangent to the geodesic leaves of $w_i$. Since $\pi_i$ are mutually orthogonal,  the metric is diagonal in the chosen coordinates
$$
g=A(x,y,z)dx^2+B(x,y,z)dz^2+C(x,y,z)dz^2.
$$
Let $d=\sqrt{A\xi^2+B\eta^2+C\zeta^2}$ be the length of the vector $v$ and $t$ be the  natural parameter on a geodesic. The geodesic equations
$$
\frac{d^2x^{\alpha}}{dt^2}+\sum_{j,k}\Gamma^{\alpha}_{jk}\frac{dx^j}{dt}\frac{dx^k}{dt}=0,
$$
where $x^1=x,\ x^2=y,\ x^3=z$ and $\Gamma^{\alpha}_{jk}$ are Christoffel symbols for  the Levi-Civita connection of our metric,
give a system of partial differential equations for $A,B,C$.
Similarly, we get a set of PDEs using the fact that $v_1=-\xi\partial_x+\eta\partial_y+\zeta\partial_z$, $v_2=\xi\partial_x-\eta\partial_y+\zeta\partial_z$, and $v_3=\xi\partial_x+\eta\partial_y-\zeta\partial_z$, obtained by applying reflections at the planes $\pi_i$ to $v$, are also Jacobi vector fields, i.e. tangent to geodesics.

We use  $\frac{dx}{dt}=\frac{\xi}{d}$, $\frac{dy}{dt}=\frac{\eta}{d}$, $\frac{dz}{dt}=\frac{\zeta}{d}$ and $\frac{d^2x}{dt^2}=\frac{1}{d}v\left(\frac{\xi}{d}\right)$, $\frac{d^2y}{dt^2}=\frac{1}{d}v\left(\frac{\eta}{d}\right)$, and $\frac{d^2z}{dt^2}=\frac{1}{d}v\left(\frac{\zeta}{d}\right)$. One can solve thus obtained system for $A_y$, $A_z$, $B_x$, $B_z$, $C_x$, $C_z$. This system simplifies for $f,j,h$, defined via
$$
A = \frac{f}{l - m},\ \ \ \ B= \frac{j}{m - k}, \ \ \ \ C = \frac{h}{k - l},
$$
namely, we have
$$
f_y=\frac{ (m - k)fj_y+l'f(f-h)}{(l - m)f + 2(m-k)j + (k-l)h}, \ \ \
f_z=\frac{(k-l)fh_z +m'f(j-f) }{(l - m)f +(m-k)j +2 (k - l)h},
$$
$$
j_x=\frac{(l - m)jf_x +k' j(h - j)}{2(l - m)f + (m-k)j +(k-l)h}, \ \ \
j_z=\frac{(k-l)jh_z+m'j(j-f)}{(l-m)f+(m-k)j  + 2(k-l)h},
$$
$$
h_x=\frac{(l-m)hf_x+k'h(h-j)}{2(l-m)f +(m-k)j +(k-l)h}, \ \ \
h_y=\frac{(m-k)hj_y +l'h(f - h)}{(l-m)f+2(m-k)j+(k-l)h}.
$$
Observe that, for definite metrics, none of the denominators can vanish  as they  are equal to length squares of some non-zero vectors. For example,  $(l - m)f + 2(m-k)j + (k-l)h=A\xi^2+B(\sqrt{2}\eta)^2+C\zeta^2$.

 One checks that
$$
\partial_z(\ln(f)-\ln(j))=\partial_y(\ln(f)-\ln(h))=\partial_x(\ln(j)-\ln(h))=0.
$$
Therefore
$$
\ln(j)=\ln(f)+\alpha (x,y),\ \ \ \ln(h)=\ln (f)+\beta (x,z), \ \ \mbox{where} \ \ \alpha(x,y)-\beta(x,z)=\gamma(y,z).
$$
From the first prolongation of the system for $f,j,h$, one obtains also
$$
\partial_x\partial_y(\ln(f)-\ln(j))=\partial_x\partial_z(\ln(f)-\ln(h))=\partial_y\partial_z(\ln(j)-\ln(h))=0.
$$
Therefore
$$
\ln(j)=\ln(f)+\bar{b}(y)-\bar{a}(x), \ \ \ \ln(h)=\ln(f)+\bar{c}(z)-\bar{a}(x)
$$
and with $a=e^{\bar{a}}$, $b=e^{\bar{b}}$, $c=e^{\bar{c}}$, $f=aS(x,y,z)$ we get  $j=bS(x,y,z)$, $h=cS(x,y,z)$.
The above system of $PDE$ for $f,j,h$ now rewrites as a completely integrable system for $S$. Moreover, holds
$$
\frac{dS}{S}=\frac{d(a(l- m) + b(m - k) + c(k - l))}{a(l- m) + b(m - k) + c(k - l)}.
$$
Hence $S=const\cdot (a(l- m) + b(m - k) + c(k - l))$ and the metric has the Stäckel form (\ref{stform}).
\hfill $\Box$\\

\section{On quadratic integrals with not normal principal congruences}\label{secprincip}
We choose a basis $\omega_1$, $\omega_2$, $\omega_3$ for $T^*_mM$ so that all $\omega_i$ are smooth, the distribution $\omega_i=0$ is orthogonal to the principal direction $n_i$ of one of the integrals $I$, and $g=\omega_1^2+\omega_2^2+\omega_3^2$. Then
$I=F\omega_1^2+G\omega_2^2+H\omega_3^2$.
If all the distributions $\omega_i=0$ are not integrable then the coefficients $K,L,M$ in the structure equations
\begin{equation}
\begin{array}{l}\label{structure}
         d\omega_1= a \omega_1 \wedge\omega_2 + b \omega_1 \wedge \omega_3 + K \omega_2 \wedge\omega_3,\\

         d\omega_2 = p \omega_2 \wedge \omega_3 + q \omega_2 \wedge \omega_1 + L \omega_3 \wedge \omega_1,\\

         d\omega_3= r \omega_3 \wedge \omega_1 + s \omega_3 \wedge \omega_2 + M \omega_1 \wedge \omega_2,
\end{array}
\end{equation}
do not vanish.

Symplectic structure on $M$ is defined by the form
$$
\Omega=d(u\omega_1+v\omega_2+w\omega_3),
$$ where $u\omega_1+v\omega_2+w\omega_3\in T^*M$.  The Hamiltonian vector field, corresponding to $\frac{1}{2}(u^2+v^2+w^2),$ is
\begin{equation}\label{V_g}
V_g=-ue_1-ve_2-we_3+\alpha \partial_u+\beta \partial_v+\gamma \partial_w,
\end{equation}
where $(e_1,e_2,e_3)$ is the basis on $T_mM$ dual to $\omega_1, \omega_2, \omega_3$ (i.e. $\omega_i(e_k)=\delta_{ik}$) and
$$
\begin{array}{l}
\alpha =(L-M)vw - u(av+bw) + qv^2 + rw^2,\\
\beta =(M-K)wu-v(pw+qu) +sw^2+ au^2,\\
\gamma=(K-L)uv-w(ru+sv) + bu^2 + pv^2.
\end{array}
$$
{\bf Proof of Theorem \ref{normal}:} Splitting the equation $V_g(Fu^2+Gv^2+Hw^2)=0$, polynomial in $u,v,w,$  with respect to these variables, one obtains a system of equations relating $F,G,H,F_i,G_i,H_i$. These eqauations allow to find all derivatives:
$$
\begin{array}{c}
F_1=G_2=H_3=0,\\
\\
F_2=2a(G-F),\ \ \ F_3=2b(H-F),\\
\\
G_1=2q(F-G), \ \ \ G_3=2p(H-G)\\
\\
H_1=2r(F-H),\ \ \ H_2=2s(G-H).
\end{array}
$$
Moreover, there is an equation not involving the derivatives:
$$
(L - M)F + (M-K)G + (K - L)H=0.
$$
Since there are 3 independent  quadratic integrals, this last equation must be trivial:
$$
M=L=K.
$$
If  $M=L=K=0$ then all the distributions $\omega_i=0$ are integrable. Moreover, we claim that $F,G,H$ are pairwise distinct. In fact, if $G\equiv F$ then the space of quadratic integrals is at most 2-dimensional.  Thus, the metric is of Stäckel type due to the result  of Di Pirro \cite{P-96}.

Consider the case  $M=L=K\ne 0$.

The form $F_1\omega_1+F_2\omega_2+F_3\omega_3$ is closed, this condition splits into 3 equations. The coefficient of $\omega_3\wedge \omega_1$ gives
$$
(2b_1-2aK - 6br)F + 2aKG + (6br - 2b_1)H=0.
$$
Since no relation between $F,G,H$ is possible and $K\ne 0$, one gets $a=0$.
Simlarly, the coefficient of $\omega_1\wedge \omega_2$ gives
$$
(6aq-2bK - 2a_1)F + (2a_1-6aq)G + 2bKH=0.
$$
Hence $b=0$. Inspecting the forms defining $dG$ and $dH$, one obtains also $p=q=r=s=0$. With vanishing $a,b,p,q,r,s$, the forms for $dF,dG,dH$ also vanish, so $F,G,H$ are constant.

Suppose that $J=Uu^2+Vv^2+Ww^2$ is another quadratic integral. Then, as we have just derived, $U,V,W$ are constant. Now the Hamiltonian vector field, corresponding to $I=\frac{1}{2}(Fu^2+Gv^2+Hw^2)$ is
$$
V_I=-uFe_1-vGe_2-wHe_3+vwK(H-G) \partial_u+wuK(F-H) \partial_v+uvK(G-F) \partial_w,
$$
and $V_I(J)=2uvwK((H-G )U + (F - H)V + (G-F)W)=0$. Hence $(H-G )U + (F - H)V + (G-F)W=0$ and the integrals
$I,J$ and $u^2+v^2+w^2$ are linearly dependent. \hfill $\Box$

  \medskip

\noindent {\bf Remark 1.} In the above proof,   the structure equations with $M=L=K$ give $K_1=K_2=K_3=0$, since $d\circ d(\omega_i)=0$. Thus the space of quadratic integrals, simultaneously diagonalizable, is three-dimensional if the structure equations are those for the  Maurer-Cartan forms of the group $SO(3)$ or $SO(2,1)$, though one can choose only two independent commuting integrals.

\medskip

\noindent {\bf Remark 2.} Continuing the study of the compatibility conditions for the case $M=L=K=0$ we easily get a system in involution for the Stäckel metric. The differentials of $F,G,H$ are written as differential forms in terms of $F,G,H$ and the coefficients of the structure equations and their derivatives. These forms must be closed. Moreover, the identities $d\circ d(\omega_1)=0$, $d\circ d(\omega_2)=0$, and $d\circ d(\omega_3)=0$ give more equations for the derivatives of the coefficients. Resolving thus obtained equations we get the following system for the structure equations:
\begin{equation}\label{Stackstructure}
\begin{array}{c}
da= 3aq\omega_1 + a_2\omega_2 + (3ap + 2bs-2ab)\omega_3,\ \ \
db= 3br\omega_1 + (3bs+2ap-2ab)\omega_2 + b_3w3,\\
\\
dp= (3pr+2bq - 2pq)\omega_1 + 3ps\omega_2 + p_3\omega_3, \ \ \
dq= q_1\omega_1 + 3aq\omega_2 + (3bq+2pr - 2pq)\omega_3,\\
\\
dr= r_1\omega_1 + (3ar + 2qs - 2rs)\omega_2 + 3rb\omega_3,\ \ \
ds= (3qs+2ar - 2rs)\omega_1 + s_2\omega_2 + 3sp\omega_3.
\end{array}
\end{equation}
These equations will be of use in what follows.

\section{Integrability by two linear integrals or by one linear and one quadratic integrals}
We call a polynomial integral {\it non-singular} if the components of the corresponding Killing tensor do not vanish simultaneously. Suppose that geodesic flow admits 2 independent non-singular integrals, linear in momenta. By Noether Theorem, this is equivalent to existence of 2 infinitesimal isometries, called also Killing vector fields. Namely, if $\sum_{i=1}^3a^i(x)p_i$ is a linear integral then $\sum_{i=1}^3a^i(x)\partial_{x^i}$ is the corresponding isometry. Let us give a criterion for existence of two commuting non-singular linear integrals in terms of geodesic flow.
\begin{lemma}\label{orto} Let $I_1=\sum_{i=1}^3a^i(x)p_i$ be a non-singular linear integral of geodesic flow in three-dimensional Riemannian manifold. Then the plane of zero directions of the integral is orthogonal to the Killing vector field.
\end{lemma}
{\it Proof:} The zero directions annihilate the form $\sum_{i,j=1}^ng_{ij}a^idx^j$, hence the claim.
\hfill $\Box$\\

Since the integrals are non-singular and commuting, the corresponding Killing vector fields can be brought in form of translations  $\partial_x,\partial_y$ in some local coordinates $x,y,z$, where $\partial_z$ is the common zero direction of integrals, orthogonal to $\partial_x,\partial_y$ due to Lemma \ref{orto}. Then, changing the coordinate $z$ if necessary, we bring the metric in the form
\begin{equation}\label{metric2I1}
g=g_{1,1}(z)dx^2+2g_{1,2}(z)dxdy+g_{2,2}(z)dy^2+dz^2.
\end{equation}

\begin{theorem}Geodesic flow on a three-dimensional Riemannian manifold $M$ admits two commuting non-singular independent integrals linear in momenta if and only if $M$ carries a normal geodesic congruence such that the restriction of the metric on each surface orthogonal to the congruence is flat.
\end{theorem}
{\it Proof:} The metric (\ref{metric2I1}) is obviously flat on surfaces $z=const$. Therefore it is enough to prove that the given condition is sufficient. The surfaces orthogonal to the congruence can be chosen as coordinate surfaces $z=const$. Moreover,  one can choose $z$ so that the unitary vector field along the geodesics is $\partial_z$ (see \cite{E-49} p. 57). Then  we can represent $M$ locally as a trivial vector bundle $M\to \mathbb{R}\times V^2$, $m\mapsto (z,v)$ where $V^2$ is a two-dimensional real vector space. Moreover, each fiber is equipped with a scalar product smoothly depending on $z$. Choosing a basis smoothly depending on $z$ we get the metric in a form (\ref{metric2I1}) with two Killing vector fields $\partial_x,\partial_y$.
\hfill $\Box$\\

If the geodesic flow is integrable by one linear and one quadratic integral then we can construct a geodesic 2-web whose leaf directions are the common zero directions of the integrals. Let us choose a basis $(\omega^1,\omega^2,\omega^3)$ in $T^*M$ so that 1) the distributions $\omega^1=0$ and $\omega^2=0$ are symmetry planes for the web directions at each point $m\in M$, 2)  the metric is $g=\omega_1^2+\omega_2^2+\omega_3^2$. In what follows we call the distributions $\omega^1=0$ and $\omega^2=0$ the {\it symmetry distributions} of the 2-web.  Then there are also here webs of maximal rank intrinsically related to the integrals.
\begin{theorem} Suppose that the geodesic flow on a three-dimensional Riemannian manifold $M$ admits two non-singular commuting integrals, one   linear in momenta and the other quadratic such that the Hamiltonian, quadratic integral, and the square of the linear integral are linearly independent and share the same principal directions.  Then the  symmetric distributions are integrable and  two 3-webs formed by symmetry distributions and one foliation of the geodesic 2-web are of maximal rank 1.
\end{theorem}
{\it Proof:} If we choose the basis $\omega^1,\omega^2,\omega^3$ as above (i.e. $\omega^1=0$ and $\omega^2=0$ are symmetry distributions and the metric is $g=\omega_1^2+\omega_2^2+\omega_3^2$), then, by Lemma \ref{orto}, the linear integral is $Rw$ for some function $R$,  the quadratic one being $Fu^2+Gv^2+Hw^2$, where the notation of Section \ref{secprincip} is used. The Hamiltonian vector field, corresponding to $Rw$ is
$$
V_1=-Re_3+(LRv - Rbu + Rrw + R_1w)\partial_u+(Rsw-KRu - Rpv  +R_2w)\partial_v+R_3w\partial_w,
$$
where $b,r,s,p,K,L$ are defined by the structure equations (\ref{structure}). Splitting the equation $V_1(u^2+v^2+w^2)=0$, polynomial in $u,v,w,$  with respect to these variables, one obtains a system of equations that allows to find all derivatives of $R$:
$$
R_1=-rR,\ \ \ R_2=-sR,\ \ \ R_3=0.
$$
Moreover, these equations imply $b=p=0$ and $R(K - L)=0$. Note that $R\ne 0$, therefore $L=K$.

Similarly, splitting the equation $V_1(Fu^2+Gv^2+Hw^2)=0$, one gets
$$
F_3=G_3=H_3=0, \ \ \ KR(F-G)=0.
$$
Since $R\ne 0$, one concludes that either $G=F$ or $K=0.$

If $G=F$ then the equation  $V_g(Fu^2+Gv^2+Hw^2)=0$, where $V_g$ is given by (\ref{V_g}) gives
$$
F_1=F_2=0,\ \ \  H_1=2r(F - H), \ \ \ H_2=2s(F - H).
$$
Thus all derivatives of $F$ and $H$ are fixed and the space of quadratic integrals, including the square of the linear one, is at most 2-dimensional.

Therefore $K=L=0$ and the symmetry distributions are integrable.  Splitting the equation  $V_g(Fu^2+Gv^2+Hw^2)=0$ one gets
$$
F_1=G_2=0,\ \ \ F_2=2a(G-F), \ \ \ G_1=2q(F-G), \ \ \ H_1=2r(F-H), \ \ \ H_2=2s(G-H),
$$
and $M(F-G)=0$. As we have just shown, holds $G\ne F$ therefore also $M=0$ and the distribution $w_3=0$ is also integrable.

Now $d\circ d(\omega_1)=0$, $d\circ d(\omega_2)=0$, and $d\circ d(\omega_3)=0$ give $a_3=q_3=0$ and $s_1=r_2 -ar + qs$ respectively. Similarly, $d\circ d(R)=0$ implies $r_3=s_3=0$, $d\circ d(F)=0$ and $d\circ d(G)=0$  imply $a_1=q_2=3aq$. Finally $d\circ d(H)=0$ gives $r_2= 3ar + 2qs - 2rs$. Thus we get an exterior differential system in involution
\begin{equation}\label{invol12}
\begin{array}{c}
d(\omega_1)= a\omega_1\wedge \omega_2,\ \ \ d(\omega_2)= q\omega_2\wedge \omega_1, \ \ \ d(\omega_3)= r\omega_3\wedge \omega_1+s\omega_3\wedge \omega_2,\\
\\
da= 3aq\omega_1 + a_2\omega_2,\ \ \ dq=q_1\omega_1 +3aq\omega_2,\\
\\
dr= r_1\omega_1 + (3ar + 2qs - 2rs)\omega_2, \ \ \ ds= (3qs+2ar - 2rs)\omega_1 + s_2\omega_2,\\
\\
dR= -rR\omega_1 - sR\omega_2,\ \ \ dF = 2a(G-F)\omega_2, \ \ \ dG= 2q(F- G)\omega_1,\\
\\
dH= 2r(F - H)\omega_1 + 2s(G - H)\omega_2.
\end{array}
\end{equation}
The 2-web defined by the common zero directions of the integrals is determined by $\omega_3=0$ and $F(\omega_1)^2+G(\omega_2)^2=0$. For real web we can assume $F>0$, $G<0$. Then in the notation of of Section \ref{secmaxrank} we take $\bar{\omega}_1=\sqrt{F}\omega_1$ and $\bar{\omega}_2=\sqrt{-G}\omega_2$. Applying Lemma \ref{2s1c} we conclude that the curvature web exists and vanishes due to equations (\ref{invol12}).
\hfill $\Box$\\

The integrability of the type described in the previous Theorem translates as well in geometric terms of webs.

\begin{theorem}
Suppose that a three-dimensional Riemannian manifold $M$ carries two  following webs:  a geodesic 4-web $w_4$ such that at each point $m\in M$ there are 3 mutually orthogonal planes $\pi_i$, $i=1,2,3$ such that
\begin{enumerate}
\item the distributions $\pi_i$ are integrable,
\item the  reflections at $\pi_i$ map one direction of the geodesic web $w_4$ to the other three,
\item the rank of  the 4-web $w_i$ formed by geodesics of one foliation of $w_4$ and by three foliations by integral manifolds of $\pi_1,\pi_2,\pi_3$ is maximal and equals to 2,
\end{enumerate}
and a smooth one-parameter family of geodesic 2-webs such that at each point $m\in M$ the symmetry distributions of each 2-web of the family coincide with two of the planes $\pi_1,\pi_2,\pi_3$.

Then the geodesic flow on  $M$ admits two commuting non-singular integrals, one   linear in momenta and the other quadratic such that the Hamiltonian, quadratic integral, and the square of the linear integral are linearly independent and share the same principal directions.
\end{theorem}
{\it Proof:} We choose again a basis $\omega_1$, $\omega_2$, $\omega_3$ for $T^*_mM$ so that  $\omega_i=0$ gives distribution $\pi_i$ and  $g=\omega_1^2+\omega_2^2+\omega_3^2$. Then by Theorem \ref{webstack} and Remark 2 of Section \ref{secprincip}, the coefficients of the structure equations (\ref{structure}) satisfy equations (\ref{Stackstructure}) with $K=L=M=0$.

Due to the permutation symmetry of these equations we can suppose that the symmetry distributions of each 2-web are $\pi_2$ and $\pi_3$. Then the vector fields tangent to the leaves of the web  are  $\xi_{\pm}=xe_1\pm ye_2$ with smooth functions $x,y$ normalized to $x^2+y^2=1$.

The connection forms $\omega^k_i$ of the Levi-Cevita connection  in the basis $e_1,e_2,e_3$, defined by
$$
\nabla e_i=\sum_k \omega^k_ie_k,
$$
are as follows:
$$
\omega^1_2=-\omega^2_1=-a\omega_1+q\omega_2,\ \ \ \omega^2_3=-\omega^3_2=-p\omega_2+s\omega_3,\ \ \
\omega^3_1=-\omega^1_3=b\omega_1-r\omega_3.
$$
Equations $\nabla_{\xi_+}(\xi_+)=\nabla_{\xi_-}(\xi_-)=0$ give for $Q:=x^2$ the following relations:
$$
Q_1= 2q(Q - 1), \ \ \ Q_2 = 2aQ,\ \ \ (p - b)Q = p.
$$
Being one-parametric,  $Q$ cannot be fixed by the last equation. Therefore $b=p=0$ and the structure equations reduce to (\ref{invol12}). Hence the claim.
\hfill $\Box$\\

  \medskip

\noindent {\bf Remark.} One cannot relax the hypothesis of the above Theorem to existence of just one 2-web instead of a family. Though existence of two such webs with the same symmetry distributions is enough to deduce that $p=b=0.$

\section{Examples}

We illustrate that the hypothesis of the presented theorems are not  technical and the general picture is more complicated. We restrict ourself  to the case of flat metric $g=dx^2+dy^2+dz^2$. The linear  integrals of geodesic flow are $p,q,r$ and $a=qz-ry,b=rx-pz,c=py-qx$, where $p,q,r$ are momenta conjugated to $x,y,z$ respectively. It is well known that  general linear and quadratic integrals are  forms in these basic integrals of degree one and two respectively. Observe that these basic integrals are Pl\"ucker coordinates of a line in space, thus verifying the equation of the Pl\"ucker quadric
$$
ap+bq+cr=0.
$$

\subsection{Linear and quadratic integral} By Euclidean motion, we can bring the isometry, corresponding to a linear integral, either to translation, or to rotation, or to screw symmetry.

\subsubsection{Translation $\partial _z$} The linear integral $I_1=r$ corresponds to this translation. The general form of a commuting  quadratic integral $I_2$ is
$$
I_2=Q(p,q,r)+Ac^2+Bcr +Ccp + Dcq,
$$
where $A,B,C,D$ are constants and $Q$ is a quadratic form in $p,q,r$.

For generic $I_2$, the subspace of quadratic integrals, spanned by $H,I_2,r^2$ and  written in velocities, does not possess a basis diagonalizing simultaneously all the integrals of the subspace.

The equations $r=I_2=0$ give (for generic $I_2$) a congruence of lines of second class: through a generic point pass 2 lines of the congruence. In fact, the lines belong to the intersection of a linear and of a quadratic complexes. One easily checks that the lines are "horizontal" (i.e. $z=const$ for each line) and that they envelope a "vertical" cylinder over a conic. By rotations around $z$-axis and translations one can simplify a bit both the congruence and the integral $I_2$.

\subsubsection{Rotation $y\partial _x-x\partial_y$} This rotation gives  linear integral $I_1=c$. The general form of a commuting  quadratic integral $I_2$ is
$$
I_2= A_1(p^2 + q^2) + A_2r^2 + B_1(a^2 + b^2) + B_2c^2 + C_1(aq - bp) + C_2cr,
$$
where $A_i,B_j,C_k$ are constants.

Again, for generic $I_2$, the subspace of quadratic integrals, spanned by $H,I_2,c^2$ and  written in velocities, does not possess a basis diagonalizing simultaneously all the integrals of the subspace.

For generic $I_2$, the equations $c=I_2=0$ give  a congruence of lines of second class. One easily checks that the congruence is symmetric by rotations,  lines are tangent to a quadric, which quadric is also symmetric by rotations, and the lines intersect the $z$-axis.

\subsubsection{Screw symmetry  $\alpha \partial_z+y\partial _x-x\partial_y$} Linear integral, corresponding to this symmetry, is $I_1=\alpha r+c$, $\alpha=const\ne 0$. The general form of a commuting  quadratic integral $I_2$ is
$$
I_2= A(p^2 + q^2) + Br^2 +  Cc^2 + Dcr,
$$
where $A,B,C,D$ are constants.

Here the geometry is even more different from the studied in this paper:  the distribution  $I_1=0$, orthogonal to the symmetry, is not integrable. In fact, the distribution is the kernel of $\omega=\alpha dz+ydx-xdy$ and $d\omega\wedge \omega=2\alpha dy\wedge dx \wedge dz\ne 0$.

For generic $I_2$, the equations $I_1=I_2=0$ give  a congruence of lines of second class. The congruence is invariant by our screw symmetry, the lines are tangent to a  cylinder of revolution around  $z$-axis.

\subsection{Two quadratic integrals}

Not any quadratic integral can be included in a 3-dimensional subspace of commuting quadratic integrals. For example, $a^2 + bp$ commutes only with integrals of 2-dimensional subspace spanned by $H$ and itself.

Consider a generic quadric
$$
\frac{x^2}{A}+\frac{y^2}{B}+\frac{z^2}{C}=1.
$$ The quadratic complex of its tangent lines, which writes as 
$$
Aa^2 +Bb^2 + Cc^2 -BCp^2-ACq^2-ABr^2=0,
$$ 
(see \cite{D-12}, p. 616) gives the quadratic  integral
$$
I_2=Aa^2 +Bb^2 + Cc^2 -BCp^2-ACq^2-ABr^2.
$$
For a member of the confocal family of quadrics
\begin{equation}\label{confocal}
\frac{x^2}{A+t}+\frac{y^2}{B+t}+\frac{z^2+t}{C+t}=1,
\end{equation}
the quadratic complex of tangent lines is $I_2+tJ_2+2Ht^2=0$, where
$$
J_2=(B + C)p^2 + (A + C)q^2 + (A + B)r^2-a^2 - b^2 - c^2.
$$ Then the subspace of quadratic integrals commuting with $I_2$ is spanned by $H,I_2,J_2$.  These integrals play the crucial role in integrating the geodesic flow on ellipsoids. (See  \cite{J-84} for classical and \cite{M-78} for modern treatment.)
The integrals of the space are diagonal in the confocal coordinates fixed by the family (\ref{confocal}).

For complete classification of Stäckel coordinate system in 3D Eucledean space see \cite{K-86,HMS-05}  and \cite{VM-22} for the Minkowski counterpart.

\section{Concluding remarks}
\subsection{Geodesically equivalent metrics}
Two metrics on a manifold $M$ are called geodesically equivalent if they have the same unparameterized geodesics. Metrics, for which there exist a nontrivial geodesically equivalent one, form an  important subclass of Stäckel metrics (see, for instance,  \cite{S-02}). They are obtained from the general formula (\ref{staeckelmetric}) by specifying a Vandermonde type matrix $\varphi_{ij}(x_i)=(\psi_i(x^i))^{j-1}$. If, for a given metric $g$, there exists a geodesically equivalent (and strictly nonproportional) metric $\bar{g}$ then the set of $n-1$ additional commuting quadratic independent integrals for geodesic flow of $g$ are obtained in a pure algebraic way from $\bar{g}$ (see \cite{MT-98}). This is a particular case of the more general result of Benenti \cite{B-92}.
\subsection{Stäckel metrics and integrable 3D billiards} Suppose that one of the distributions $\pi_i$,  described in Proposition \ref{mirrors},  has a closed integral surface $\mathcal{B}$.  Consider a point $m$ on this surface and one of the 4 geodesics $L_k$ whose directions are  defined by equations (\ref{webdirectioneq}).  If a ball arrives at $m$ along $L_k$ then,  reflected elastically from the wall $\mathcal{B}$,  it continues to move along the geodesic $L_j$, whose direction is obtained by reflecting $L_k$ in $\pi_i$.  By Proposition \ref{mirrors},  the direction of $L_j$ at $m$ is defined by the same equations (\ref{webdirectioneq}).  For definite metrics $g$ and any direction $\tau$, we can choose  the constants $\lambda,\mu$ so that (\ref{webdirectioneq}) are satisfied. Therefore we can  send the ball in any given direction,  playing the 3D billiard with the wall $\mathcal{B}$.  Moreover,  this billiard is integrable. In fact,  the values of $\lambda,\mu$   remains constant after reflection.

The geometric properties of the constructed billiard are analogous to the classical 3D billiard in quadrics. For instance, for  $c_2=\frac{M(z_0)}{G(z_0)}$, $c_3=\frac{R(z_0)}{G(z_0)}$ the zero directions of the integral $g+c_2I_2+c_3I_3$ are tangent to the coordinate surface $z=z_0$. Thus,  the Stäckel coordinates play the role of elliptic coordinates for billuards in quadrics (for more details on billiards in quadrics see \cite{DR-11}). Equations (\ref{directions}), defining the zero directions of pair of integrals, are analogous of equations written by Jacoby in terms of hyperelliptic funcions in elliptic coordinates (see \cite {J-84} or \cite{DR-11} p.168).

The presented construction of integrable billiards in 3D via integrable geodesic flow generalizes the one described in \cite{Aq-21} for two-dimensional case.

\section*{Acknowledgments}
The author thanks Boris Kruglikov for useful discussions. This research was supported by FAPESP grant \# 2022/12813-5 and partially by Troms\o\  Research Foundation.   The author also thanks the personnel of the Department of Mathematics and Statistics of UiT the Arctic University, where this study was initiated, for their warm hospitality.

\end{document}